\documentclass[letterpaper, 10 pt, conference]{ieeeconf}  % Comment this line
\usepackage{xcolor}
\definecolor{commentgreen}{RGB}{0, 153, 0}

\usepackage{amsthm}
\usepackage{amsmath}               
{
  \theoremstyle{plain}
  \newtheorem{assumption}{Assumption}
  \newtheorem{definition}{Definition}
  \newtheorem{theorem}{Theorem}
  \newtheorem{corollary}{Corollary}
  \newtheorem{lemma}{Lemma}
  
}
\usepackage{amssymb}
\usepackage{mathtools}
\usepackage{graphicx}
\usepackage{tabularx}
\usepackage[]{hyperref}
\usepackage{xcolor}
\usepackage{algorithm, algorithmicx}
\usepackage[noend]{algpseudocode}
\usepackage{pgffor}
\usepackage{etoolbox}
\usepackage{multirow}
\usepackage{ifthen}
\usepackage{subcaption}

\algrenewcommand\algorithmicindent{0.9em}%
\algnewcommand{\algorithmicgoto}{\textbf{go to} line}%
\algnewcommand{\Goto}[1]{\algorithmicgoto~\ref{#1}}%

\newcommand{\BEAS}{\begin{eqnarray*}}
\newcommand{\EEAS}{\end{eqnarray*}}
\newcommand{\BEA}{\begin{eqnarray}}
\newcommand{\EEA}{\end{eqnarray}}
\newcommand{\BEQ}{\begin{equation}}
\newcommand{\EEQ}{\end{equation}}
\newcommand{\BIT}{\begin{itemize}}
\newcommand{\EIT}{\end{itemize}}
\newcommand{\BNUM}{\begin{enumerate}}
\newcommand{\ENUM}{\end{enumerate}}

\newcommand{\BA}{\begin{array}}
\newcommand{\EA}{\end{array}}

\newcommand{\reals}{\mathbb R}
\newcommand{\ball}{\mathbb B}

\newcommand{\binary}{\mathbb I}

\newcommand{\Spec}{\mathrm{spec}}
\newcommand{\Co}{{\mathop \mathrm{co}}}

{\begin{quote}}{\end{quote}}

\newcounter{oursection}

\DeclareMathOperator*{\find}{find}

\makeatletter
\def\optimize{\@ifnextchar[{\@with}{\@without}}
\def\@with[#1]#2#3#4#5{
  \ifthenelse{\equal{#2}{}}{
    \begin{aligned}
      #3_{#4}\text{ } & #5 \\
      \mathrm{s.t.}\text{ } & #1
    \end{aligned}
  }{
    \begin{aligned}
      #2 & = & #3_{#4}\text{ } & #5 \\
         &   & \mathrm{s.t.}\text{ } & #1
    \end{aligned}
  }
}
\def\@without#1#2#3#4{
  \ifthenelse{\equal{#1}{}}{
    \begin{aligned}
      #2_{#3}\text{ } & #4
    \end{aligned}
  }{
    \begin{aligned}
      #1 & = & #2_{#3}\text{ } & #4
    \end{aligned}
  }
}
\makeatother

\definecolor{darkolivegreen}{rgb}{0.33, 0.42, 0.18}

\newcommand{\Behcet}{Beh\c{c}et}
\newcommand{\Acikmese}{A\c{c}{\i}kme\c{s}e}

\newcommand{\NPcomplete}{$\mathcal{NP}$-complete}

\newcommand{\transp}{{\scriptscriptstyle\mathsf{T}}}

\IEEEoverridecommandlockouts                              % This command is only needed if 
\usepackage{graphics} % for pdf, bitmapped graphics files

\title{\LARGE \bf Partition-based Feasible Integer Solution Pre-computation for
  Hybrid Model Predictive Control}

\author{Danylo Malyuta$^{1}$, \Behcet{} \Acikmese{}$^{2}$, Martin Cacan$^{3}$
  and David S. Bayard$^{4}$
  \thanks{$^{1}$Ph.D. student, W.E. Boeing Department of Aeronautics \&
    Astronautics, University of Washington, Seattle, WA 98195, USA
    \texttt{danylo@uw.edu}}%
  \thanks{$^{2}$Professor, W.E. Boeing Department of Aeronautics \&
    Astronautics, University of Washington, Seattle, WA 98195, USA
    \texttt{behcet@uw.edu}}%
  \thanks{$^{3}$Guidance and Control Analyst, Jet Propulsion Laboratory,
    California Institute of Technology, Pasadena, CA 91109, USA
    \texttt{martin.cacan@jpl.nasa.gov}}%
  \thanks{$^{4}$Senior Research Scientist, Jet Propulsion Laboratory, California
    Institute of Technology, Pasadena, CA 91109, USA
    \texttt{david.s.bayard@jpl.nasa.gov}}%
}

\begin{document}

\maketitle
\thispagestyle{empty}
\pagestyle{empty}

%%%%%%%%%%%%%%%%%%%%%%%%%%%%%%%%%%%%%%%%%%%%%%%%%%%%%%%%%%%%%%%%%%%%%%%%%%%%%%%%
\begin{abstract}
  For multiparametric mixed-integer convex programming problems such as those
  encountered in hybrid model predictive control, we propose an algorithm for
  generating a feasible partition of a subset of the parameter space. The result
  is a static map from the current parameter to a suboptimal integer solution
  such that the remaining convex program is feasible. Convergence is proven with
  a new insight that the overlap among the feasible parameter sets of each
  integer solution governs the partition complexity. The partition is stored as
  a tree which makes querying the feasible solution efficient. The algorithm can
  be used to warm start a mixed integer solver with a real-time guarantee or to
  provide a reference integer solution in several suboptimal MPC schemes. The
  algorithm is tested on randomly generated systems with up to six states,
  demonstrating the effectiveness of the approach.
\end{abstract}
%%%%%%%%%%%%%%%%%%%%%%%%%%%%%%%%%%%%%%%%%%%%%%%%%%%%%%%%%%%%%%%%%%%%%%%%%%%%%%%%

%%%%%%%%%%%%%%%%%%%%%%%%%%%%%%%%%%%%%%%%%%%%%%%%%%%%%%%%%%%%%%%%%%%%%%%%%%%%%%%%
\section{Introduction}
\label{introduction}

Model predictive control (MPC) is a discrete-time control technique in which a
receding horizon optimization problem is solved in order to determine the
optimal control input at each time step. Advanced formulations of MPC include
hybrid MPC (HMPC) and robust MPC (RMPC)
\cite{Mayne2000,Bemporad2007,Mayne2014}. HMPC handles systems like chemical
powerplants, pipelines and aerospace vehicles whose dynamics involve either
explicit discrete switches such as valves
\cite{Bemporad1999a,OcampoMartinez2007} or have nonlinearities that can be
appropriately modeled via a piecewise affine approximation
\cite{Blackmore2012,Schouwenaars2006}. RMPC handles systems that are affected by
uncertainties such as in their dynamics, in the state estimate and in the input
\cite{Bemporad2007,Malyuta2019}. Many practical applications call for a combined
control of uncertain hybrid systems which requires solving a convex
mixed-integer nonlinear program (MINLP)
\cite{Richards2003,Hen2002,Corona2006}. While possible on powerful hardware,
on-line MINLP solution is both slow and \NPcomplete{} \cite{Bemporad1999a},
meaning that there is generally no real-time performance guarantee.

To improve MPC on-line computational efficiency, some authors have worked on
explicit MPC techniques which reduce on-line computational demand by
pre-computing off-line all or part of the optimal solution. On-line it typically
remains to evaluate a piecewise affine (PWA) function. Various explicit MPC
methodologies have been proposed
\cite{Alessio2009,Bemporad2006,Pistikopoulos2012}. When the MPC law is a linear
or a quadratic program, an exact explicit law can be obtained by solving a
multiparametric program. Exact solutions for more general programs are generally
not possible due to non-convexity of common active constraint regions
\cite{Bemporad2006b}. Instead, approximate solutions have been proposed via
local linearization \cite{Pistikopoulos2007a,Oberdieck2017} or via optimal cost
bounding by PWA functions over simplices
\cite{Bemporad2006b,MunozDeLaPena2004,MunozDeLaPena2006} or hyperrectangles
\cite{Johansen2004}. An approximate explicit solution to mixed-integer quadratic
programs has been proposed based on difference of convex programming
\cite{Alessio2006} and for MINLPs based on local linearization and primal/master
subproblems \cite{Dua1998,Rowe2003}.

Multiparametric programming, however, is restricted to relatively low
dimensional systems due to the worst-case exponential partition
complexity. Several authors in hybrid MPC have therefore suggested to retain
on-line mixed-integer programming and either to reduce the integer variable's
degrees of freedom \cite{Ingimundarson2007} or to abort the solver at a
suboptimal solution when computation time becomes excessive
\cite{Bemporad1999a}. The former solution, however, has no rigorous way of
selecting a reference integer solution while the latter relies on the ability to
use the previous time step's solution to warm start the mixed-integer solver,
which is not always possible such as, for example, in some robust MPC schemes
\cite{Malyuta2019}.

To address the issue of guaranteed real-time computation of a feasible integer
solution in the general setting, our main contribution is a novel partitioning
algorithm which pre-computes feasible integer solutions in a polytopic subset of
the state space. This partition is stored as a tree which can be queried in
polynomial time. As a result, the partition provides a guaranteed real-time warm
start capability to the mixed-integer solver and thus is helpful for
\cite{Bemporad1999a,Ingimundarson2007}. Our second contribution is a convergence
proof of the algorithm which for the first time in literature considers an
overlap property as being a driver of partition complexity.

The paper is organized as follows. In Section~\ref{sec:problem_formulation} the
scope of MPC formulations that our algorithm can handle is defined as a generic
MINLP. The algorithm is then described in Section~\ref{sec:phase1} and its
convergence, complexity and use-cases are discussed in
Section~\ref{sec:properties}. The algorithm is tested on a large set of randomly
generated dynamical systems with up to 6 states and 21 integer variable degrees
of freedom, indicating that it is robust and can scale to medium dimensional
systems. Section~\ref{sec:future} outlines future research directions and is
followed by some concluding remarks in Section~\ref{sec:conclusion}.

%%% Local Variables:
%%% mode: latex
%%% TeX-master: "root"
%%% End:

\section{Problem Formulation}
\label{sec:problem_formulation}

This section defines a template MINLP that generates all MPC laws that our
algorithm can handle. Because MPC is fundamentally an optimization problem, we
do this without specific mention of a receding-horizon optimal control problem.

We use the following notation. $\reals$ denotes the set of reals,
$\binary\triangleq\{0,1\}$ the binary set and
$\ball\triangleq\{x:\|x\|_2\le 1\}$ the unit ball. Unless otherwise specified,
matrices are uppercase (e.g. $A$), scalars, vectors and functions are lowercase
(e.g. $x$) and sets are calligraphic uppercase (e.g. $\mathcal S$). We use
$1_n\in\reals^n$ to denote the vector of ones and $I_n\in\reals^{n\times n}$ to
denote the identity matrix. The scalar $\ell$ is a placeholder for some integer
value. The operator
$\mathrm{diag}(\{x_i\}_{i=1}^\ell)\in\reals^{\ell\times\ell}$ creates a diagonal
matrix with value $x_i$ at row and column $i$ and zero otherwise. The constraint
$M\succ(\succeq) 0$ means that $M\in\reals^{n\times n}$ is positive
(semi)definite. $\mathcal S^{\mathrm{c}}$, $\partial\mathcal S$,
$c^{\mathcal S}$ and $\mathcal V(\mathcal S)$ denote respectively the
complement, boundary, barycenter and vertices of $\mathcal S$ (with the latter
only relevant when $\mathcal S$ is a polytope). Given
$\mathcal A\subseteq\reals^n$, $s\in\reals$ and $b\in\reals^n$,
$\mathcal A+b\triangleq\{a+b\in\reals^n:a\in\mathcal A\}$ translates and
$s\mathcal A\triangleq\{sa: a\in\mathcal A\}$ scales the set $\mathcal A$. The
shorthand $a:b$ denotes the integer sequence $a,\dots,b$.

The following multiparametric conic MINLP generates all MPC formulations that
our algorithm can handle:
\begin{equation}
  \label{eq:minlp}
  \optimize[
  g(\theta,x,\delta)=0, \\
  &&& h(\theta,x,\delta)\in\mathcal K, \\
  &&& \delta\in\binary^m,
  ]{V^*(\theta)}{\min}{x,\delta}{f(\theta,x,\delta)}
  \tag{P$_\theta$}
\end{equation}
where $\theta\in\reals^p$ is the parameter, $x\in\reals^n$ is the decision
vector and $\delta\in\binary^m$ is a binary vector called the
\textit{commutation}. The cost function
$f:\reals^p\times\reals^n\times\binary^m\to\reals$ is jointly convex in its
first two arguments while the constraint functions
$g:\reals^p\times\reals^n\times\binary^m\to\reals^{\ell}$ and
$h:\reals^p\times\reals^n\times\binary^m\to\reals^d$ are affine in their first
two arguments. The functions can be nonlinear in the last argument. The convex
cone
$\mathcal K=\mathcal C_1^{c_1}\times\cdots\times\mathcal C_{q}^{c_q}\subset\reals^d$
is a Cartesian product of $q$ convex cones where
$\mathcal C_i^{c_i}\subset\reals^{c_i}$ and $d\triangleq \sum_{i=1}^{q}c_i$,
similar to \cite{Domahidi2013}. Examples include the positive orthant
$\reals_+^n$, the second-order cone
$\mathcal Q^{\ell}=\{(t,z)\in\reals\times\reals^{\ell-1}: t\ge\|z\|_2\}$ and the
positive semidefinite cone
$\mathcal S_+^\ell=\{Z\in\reals^{\ell\times\ell} : Z\succeq 0\}$. We also define
the following fixed-commutation multiparametric conic NLP:
\begin{equation}
  \label{eq:nlp}
  \optimize[
  g(\theta,x,\delta)=0, \\
  &&& h(\theta,x,\delta)\in\mathcal K,
  ]{V^*_\delta(\theta)}{\min}{x}{f(\theta,x,\delta)}
  \tag{P$_{\theta}^\delta$}
\end{equation}
which corresponds to (\ref{eq:minlp}) where $\delta$ has been fixed,
i.e. assigned a specific value. The following definitions closely follow
\cite{Bemporad2006b}.

\begin{definition}
  The feasible parameter set $\Theta^*\subset\reals^p$ is the set of all
  $\theta$ parameters for which (\ref{eq:minlp}) is feasible.
\end{definition}

\begin{definition}
  The fixed-commutation feasible parameter set
  $\Theta_{\delta}^*\subset\reals^p$ is the set of all $\theta$ parameters for
  which (\ref{eq:nlp}) is feasible.
\end{definition}

\begin{definition}
  The feasible commutation map $f_\delta:\Theta^*\to\binary^m$ maps
  $\theta\in\Theta^*$ to a commutation $\delta$ such that
  $\theta\in\Theta^*_\delta$.
\end{definition}

This paper presents an algorithm for computing $f_\delta$ over a subset of its
domain $\Theta\subseteq\Theta^*$. We shall assume in Section~\ref{sec:phase1}
that the set $\Theta$ is available as a convex and full-dimensional polytope in
vertex representation. Section~\ref{subsec:Theta} discusses how one might obtain
$\Theta$. It is implicitly assumed throughout this paper that (\ref{eq:minlp})
and all related problems are appropriately scaled. We suggest a per-axis unit
scaling such that $\max_{\theta\in\Theta}|e_i^\transp\theta|=1$
$\forall i=1,\dots,p$ where $\{e_i\}_{i=1}^p$ is the standard Eucledian basis.

%%% Local Variables:
%%% mode: latex
%%% TeX-master: "root"
%%% End:

\section{Feasible Commutation Map Computation}
\label{sec:phase1}

This section proposes an algorithm for computing $f_\delta$. The main idea is to
generate a coarse simplicial partition
$\mathcal R=\{(\mathcal R_i,\delta_i)\}_{i=1}^{|\mathcal R|}$ such that
$\Theta=\bigcup_{i=1}^{|\mathcal R|}\mathcal R_i$ and where each cell
$\mathcal R_i$ is associated with a fixed commutation $\delta_i$ that is
feasible everywhere in $\mathcal R_i$,
i.e. $\mathcal R_i\subseteq\Theta^*_{\delta_i}$. We then define:
\begin{equation}
  \label{eq:feasible_commutation_function}
  f_\delta(\theta)=\delta_i\text{ such that }\theta\in\mathcal R_i.
\end{equation}

\begin{lemma}
  \label{lemma:convexity}
  For any fixed $\delta\in\binary^m$, $\Theta^*_\delta$ is a convex set and
  $V^*_\delta$ is a convex function.
  \begin{proof}
    Suppose $\theta',\theta''\in\Theta^*_\delta$. Let
    $\alpha',\alpha''\in [0,1]$, $\alpha'+\alpha''=1$ and
    $\theta=\alpha'\theta'+\alpha''\theta''$. Since $g,h$ are affine in their
    first argument and $\mathcal K$ is a convex cone:
    \begin{equation*}
      \begin{split}
        g(\theta,x,\delta)=
        \alpha'g(\theta',x,\delta)+\alpha''g(\theta'',x,\delta)=0, \\
        h(\theta,x,\delta)=
        \alpha'h(\theta',x,\delta)+\alpha''h(\theta'',x,\delta)\in\mathcal K,
      \end{split}
    \end{equation*}
    so $\theta\in\Theta^*_\delta$ which is thus a convex set. Next,
    (\ref{eq:nlp}) is a minimization of $f$ over a convex set in $x$ which
    preserves convexity in $\theta$ by the joint convexity property of $f$
    \cite{Boyd2004}. Thus, $V^*_\delta$ is a convex function.
  \end{proof}
\end{lemma}

\begin{algorithm}
  \centering
  \begin{algorithmic}[1]
    \State $\mathcal R\gets\emptyset$, $\bar\Theta\gets \Theta$
    \For{all $\delta\in \binary^m$}
    \State $\mathcal R\gets\{(\mathcal R',\delta):
    \mathcal R'\in\bar\Theta\cap\Theta^*_\delta\}\cup\mathcal R$
    \State $\bar\Theta\gets\bar\Theta\setminus\Theta^*_\delta$
    \If{$\bar\Theta=\emptyset$}
    \State STOP
    \EndIf
    \EndFor
    \caption{Brute force $f_\delta$ computation.}
    \label{alg:phase1bf}
  \end{algorithmic}
\end{algorithm}

Since $\Theta^*_\delta$ is convex, an arbitrarily precise inner approximation of
it can be found \cite{Dueri2016}. A conceptually trivial method for generating
$\mathcal R$ is given by Algorithm~\ref{alg:phase1bf}. The set difference and
intersection operations in Algorithm~\ref{alg:phase1bf} are element-wise
\cite{Baotic2009}. The idea is to exploit the ability to inner approximate
$\Theta^*_\delta$ to procedurally ``cover'' $\Theta$ by intersecting it with
fixed-commutation feasible parameter sets.

The filling problem is combinatorial, however, such that in the worst case all
$2^m$ possible values of $\delta$ are needed, excepting those that are
infeasible directly due to the constraints in (\ref{eq:minlp}). Furthermore,
accurate polytopic inner approximation of $\Theta^*_\delta$ in
higher-dimensional spaces than about $\reals^4$ suffers from excessive vertex
count \cite{Dueri2016}. Last but not least, the set intersection and set
difference operations used by Algorithm~\ref{alg:phase1bf} have poor numerical
properties such as creating badly conditioned (i.e. quasi lower-dimensional)
polytopes. One would like instead an algorithm which may 1) potentially avoid
exploring all $2^m$ combinations for $\delta$, 2) minimizes vertex count and 3)
uses only numerically robust operations. We thus propose
Algorithm~\ref{alg:phase1} in which the first requirement is addressed via
(\ref{eq:maxvol}) discussed below, the second by using only simplices (which
have the lowest vertex count among full-dimensional polytopes) and the third by
working solely in the vertex representation which is much more numerically
robust than the halfspace representation and avoids the numerically fragile
operation of converting between the two representations.

\begin{algorithm}
  \centering
  \begin{algorithmic}[1]
    \State Create empty tree with open leaf $\Theta$ as root
    \State $\mathcal S\gets\mathtt{delaunay}(\mathcal V(\Theta))$
    \label{alg:phase1:line:initialdelaunay}
    \State Add child open leaves $\mathcal S_i$
    $\forall\mathcal S_i\in\mathcal S$
    \label{alg:phase1:line:init}
    \While{any open leaf exists}
    \label{alg:phase1:line:whileloop}
    \State $\mathcal R\gets\text{the first open leaf}$
    \label{alg:phase1:line:getfirstopenleaf}
    \If{(\ref{eq:minlp}) is infeasible for $\theta=c^{\mathcal R}$}
    \label{alg:phase1:line:barycenterinfeascheck}
    \State STOP, $\Theta^*\setminus\Theta\ne\emptyset$
    \label{alg:phase1:line:stop}
    \Else
    \State $\hat\delta\gets\text{solve (\ref{eq:maxvol})}$
    \label{alg:phase1:line:findmaxvolcommutation}
    \If{(\ref{eq:maxvol}) is infeasible}
    \State
    $\bar v,\bar v'\gets\arg\max_{v,v'\in\mathcal V(\mathcal R)}\|v-v'\|_2$
    \label{alg:phase1:line:longestedge}
    \State $v_{\mathrm{mid}}\gets (\bar v+\bar v')/2$
    \State Add child open leaf
    $\Co\{(\mathcal V(\mathcal R)\setminus\{\bar v\})\cup\{v_{\mathrm{mid}}\}\}$
    \State Add child open leaf
    $\Co\{(\mathcal V(\mathcal R)\setminus\{\bar v'\})\cup\{v_{\mathrm{mid}}\}\}$
    \label{alg:phase1:line:splitalonglongestedge2}
    \Else
    \label{alg:phase1:line:vertexfeasibilitytest}
    \State Replace leaf with closed leaf $(\mathcal R,\hat\delta)$
    \label{alg:phase1:line:closeleaf}
    \EndIf
    \EndIf
    \EndWhile
    \caption{Proposed computation of $f_\delta$.}
    \label{alg:phase1}
  \end{algorithmic}
\end{algorithm}

Algorithm~\ref{alg:phase1} creates a simplicial partition of $\Theta$ as
follows. The partition is stored as a tree whose leaves are cells
$(\mathcal S,\delta)$ storing the set $\mathcal S$ and associated commutation
$\delta$. Non-leaf nodes in the tree store just the sets and make evaluating
(\ref{eq:feasible_commutation_function}) more efficient (see
Section~\ref{subsec:complexity}). A ``closed leaf'' refers to a cell that will
be a leaf in the final tree while an ``open leaf'' will be further partitioned
at the next iteration and thus will be merely a non-leaf node in the final
tree. In the algorithm, we refer to nodes directly by their contents,
i.e. $\mathcal S$ for open and $(\mathcal S,\delta)$ for closed leaves. On
line~\ref{alg:phase1:line:init} the tree root is initialized to $\Theta$ and,
since generally $\Theta$ is not a simplex, Delaunay triangulation is first
applied \cite[Section~9.3]{deBerg2008}. The tree is then iterated on
line~\ref{alg:phase1:line:whileloop} in a depth-first manner until no open
leaves are left. By doing a depth-first search,
Assumption~\ref{assumption:positiveoverlap} in
Section~\ref{subsec:phase1convergence} is disproved more quickly in case that it
does not hold, such that the algorithm fails with less wasted time.

An open leaf $\mathcal R$ in the tree is selected on
line~\ref{alg:phase1:line:getfirstopenleaf}. First, it is checked whether
(\ref{eq:minlp}) is feasible at its barycenter. If not, this point is a
certificate of infeasibility of (\ref{eq:minlp}) in $\Theta$, in which case
Section~\ref{subsec:Theta} should be consulted. If however (\ref{eq:minlp}) is
feasible, then
lines~\ref{alg:phase1:line:findmaxvolcommutation}-\ref{alg:phase1:line:closeleaf}
partition $\mathcal R$ into at most $2$ simplices. First, it is checked whether
$\mathcal R$ is fully contained inside some particular $\Theta^*_\delta$. The
following lemma is used for this purpose.

\begin{lemma}
  \label{lemma:containment}
  $\mathcal R\subseteq\Theta^*_{\hat\delta}$ $\Leftrightarrow$ (\ref{eq:nlp}) is
  feasible $\forall\theta\in\mathcal V(\mathcal R)$ and $\delta=\hat\delta$.
  \begin{proof}
    $(\Rightarrow)$ Since $\mathcal R\subseteq\Theta^*_{\hat\delta}$,
    $\theta\in\mathcal V(\mathcal R)\Rightarrow\theta\in\Theta_{\hat\delta}^*$. Since
    $\Theta_{\hat\delta}^*$ is the fixed-commutation feasible parameter set,
    (\ref{eq:nlp}) is by definition feasible. $(\Leftarrow)$ Any $\theta$ such
    that (\ref{eq:nlp}) is feasible satisfies, by definition,
    $\theta\in\Theta^*_{\hat\delta}$. By Lemma~\ref{lemma:convexity}, since
    $\Theta^*_{\hat\delta}$ is convex,
    $\Co\{\theta\in\mathcal V(\mathcal R)\}\equiv\mathcal R\subseteq\Theta^*_{\hat\delta}$.
  \end{proof}
\end{lemma}

Using Lemma~\ref{lemma:containment}, one can efficiently check if
$\mathcal R\subseteq\Theta^*_\delta$ for some $\delta$ via the following
feasibility MINLP:
\begin{equation}
  \label{eq:maxvol}
  \optimize[
  g(\theta,x_\theta,\delta)=0 \quad & \forall \theta\in\mathcal V(\mathcal R), \\
  &&& h(\theta,x_\theta,\delta)\in\mathcal K& \forall \theta\in\mathcal V(\mathcal R), \\
  &&& \delta\in\binary^m,
  ]{\hat \delta(\mathcal R)}{\find}{}{\delta}
  \tag{V$^{\mathcal R}$}
\end{equation}
where $x_\theta$ denotes a feasible decision vector corresponding to the
particular value of $\theta$. Problem (\ref{eq:maxvol}) can be solved in the
standard fashion as a MINLP and the found feasible commutation can be associated
with $\mathcal R$ and the leaf can be subsequently closed. If (\ref{eq:maxvol})
is infeasible, however, then $\mathcal R$ is not fully contained in any
$\Theta^*_\delta$. In this case, $\mathcal R$ is split into two smaller
simplices at the midpoint of its longest edge. As explained in
Section~\ref{subsec:phase1convergence}, this yields a volume reduction that
necessarily leads to convergence if Assumption~\ref{assumption:positiveoverlap}
holds.

%%% Local Variables:
%%% mode: latex
%%% TeX-master: "root"
%%% End:

\section{Properties}
\label{sec:properties}

\subsection{Convergence}
\label{subsec:phase1convergence}

The main result of this section is Theorem~\ref{theorem:phase1convergence} which
guarantees convergence of Algorithm~\ref{alg:phase1} under an assumption.

\begin{definition}
  \label{definition:overlap}
  Let
  $\Delta\triangleq\{\delta\in\binary^m:\Theta^*_\delta\cap\Theta\ne\emptyset\}$. The
  largest value $\kappa\in\reals_+$ such that $\forall\theta\in\Theta$
  $\exists\delta\in\Delta$ such that
  $(\kappa\ball+\theta)\setminus(\Theta^*\cap\Theta)^{\mathrm{c}}\subseteq\Theta^*_\delta$
  is called the overlap.
\end{definition}

\begin{assumption}
  \label{assumption:positiveoverlap}
  The overlap is positive, i.e. $\kappa>0$.
\end{assumption}

The overlap depends on (\ref{eq:minlp}) and the choice of
$\Theta$. Assumption~\ref{assumption:positiveoverlap} implies that a non-zero
overlap between fixed-commutation feasible parameter sets exists everywhere near
$\partial\Theta^*_\delta$ $\forall\delta\in\Delta$. This ``fuzzy'' commutation
transition property is instrumental for the convergence proof in
Theorem~\ref{theorem:phase1convergence}.

\begin{theorem}
  \label{theorem:phase1convergence}
  If Assumption~\ref{assumption:positiveoverlap} holds then
  Algorithm~\ref{alg:phase1} either converges or fails in a finite number of
  iterations.
  \begin{proof}
    Let $\mathcal R_k$ be the leaf chosen at the $k$-th call of
    line~\ref{alg:phase1:line:getfirstopenleaf}. Whenever $\mathcal R_k$ is not
    closed, it can be shown that the partition along its longest edge on
    lines~\ref{alg:phase1:line:longestedge}-\ref{alg:phase1:line:splitalonglongestedge2}
    halves the volume, therefore
    $\lim_{k\to\infty}\mathrm{vol}(\mathcal R_k)=0$. Since the longest edge
    length is also halved, $\exists k$ large enough such that
    $\mathcal R_k\subseteq(\kappa\ball+c^{\mathcal R_k})$. Two possibilities
    exist: 1)
    $\mathcal R_k\subseteq(\kappa\ball+c^{\mathcal R_k})\setminus(\Theta^*\cap\Theta)^{\mathrm{c}}$
    or 2) $\mathcal R_k\cap(\Theta^*)^{\mathrm{c}}\ne\emptyset$. In the first
    case, the $\hat\delta$ picked on
    line~\ref{alg:phase1:line:findmaxvolcommutation} is then the one feasible
    $\forall\theta\in\mathcal V(\mathcal R_k)$. Leaf $\mathcal R_k$ is therefore
    closed on line~\ref{alg:phase1:line:closeleaf}. By this logic, for a large
    enough (but finite) $k$ all regions that do not intersect the infeasible
    parameter set $(\Theta^*)^{\mathrm{c}}$ get closed. If the second case does
    not occur, the algorithm terminates.

    In the second case, recall that
    $\lim_{k\to\infty}\mathrm{vol}(\mathcal R_k)=0$. Since
    $\mathcal R_k\cap(\Theta^*)^{\mathrm{c}}\ne 0$, in a finite number of
    iterations $c^{\mathcal R_k}\notin\Theta^*$ so
    line~\ref{alg:phase1:line:barycenterinfeascheck} will evaluate to true and
    the algorithm will fail on line~\ref{alg:phase1:line:stop}.
  \end{proof}
\end{theorem}

\begin{corollary}
  \label{corollary:phase1nonconvergence}
  If Assumption~\ref{assumption:positiveoverlap} does not hold then
  Algorithm~\ref{alg:phase1} does not converge.
  \begin{proof}
    If Assumption~\ref{assumption:positiveoverlap} does not hold then there
    exists a region $\Theta'\subseteq\Theta$ such that
    $\forall\theta\in\Theta'$,
    $(\kappa\ball+\theta)\setminus(\Theta^*\cap\Theta)^{\mathrm{c}}\subseteq\Theta^*_\delta$
    for some $\delta\in\Delta\Leftrightarrow\kappa=0$. This, however, implies
    that the only simplex that would validate Lemma~\ref{lemma:containment} is a
    lower-dimensional one, i.e. with zero volume. Since this occurs at iteration
    $k=\infty$, Algorithm~\ref{alg:phase1} does not converge.
  \end{proof}
\end{corollary}

Theorem~\ref{theorem:phase1convergence} and
Corollary~\ref{corollary:phase1nonconvergence} suggest that $\kappa>0$ is not
only necessary and sufficient for convergence but that $\kappa$ also drives the
convergence rate. A small $\kappa$ implies a high iteration count $k$ until
Assumption~\ref{assumption:positiveoverlap} guarantees leaf closure. We call a
MINLP with large $\kappa$ ``well-conditioned'' and Algorithm~\ref{alg:phase1}
will converge more quickly with a rate that is derived in
Corollary~\ref{corollary:phase1cellcountcomplexity} of
Section~\ref{subsec:complexity}.

\subsection{Complexity}
\label{subsec:complexity}

In this section we analyze the complexity of Algorithm~\ref{alg:phase1} in terms
of the partition cell count as well as the on-line evaluation complexity.

\begin{theorem}
  \label{theorem:phase1treedepthcomplexity}
  The maximum tree depth $\tau$ of Algorithm~\ref{alg:phase1} is
  $\mathcal O(p^2\log(\kappa^{-1}))$.
  \begin{proof}
    Algorithm~\ref{alg:phase1} reduces search space volume by halving the
    longest edge length on
    lines~\ref{alg:phase1:line:longestedge}-\ref{alg:phase1:line:splitalonglongestedge2}. Suppose
    that $l_0$ is the longest edge length of a simplex
    $\mathcal R\subset\reals^p$, then its length is $l_k\triangleq l_0/2^k$
    after $k$ divisions. We wish to determine the number of divisions necessary
    until $\mathcal R\subseteq\kappa\ball+c^{\mathcal R}$ and gets closed by
    Theorem~\ref{theorem:phase1convergence}. This approximately requires
    $l_k\le\kappa\Rightarrow k\ge\log_2(l_0/\kappa)$. Since $\mathcal R$ has
    $p(p+1)/2$ edges then an approximate number of required subdivisions,
    i.e. the depth of the partition tree, is given by:
    \begin{equation}
      \label{eq:minkexact}
      \tau = \left\lceil\frac{p(p+1)\log_2(l_0/\kappa)}{2}\right\rceil,
    \end{equation}
    which yields $\tau=\mathcal O(p^2\log(\kappa^{-1}))$.
  \end{proof}
\end{theorem}

\begin{corollary}
  \label{corollary:phase1cellcountcomplexity}
  The maximum tree leaf count $\eta$ of Algorithm~\ref{alg:phase1} is
  $\mathcal O(2^{p^2\log(\kappa^{-1})})$.
  \begin{proof}
    In the worst case, Algorithm~\ref{alg:phase1} generates a full binary tree
    of depth $\tau=\mathcal O(p^2\log(\kappa^{-1}))$ according to
    Theorem~\ref{theorem:phase1treedepthcomplexity}, neglecting the first layer
    where the node count depends on $|\mathcal S|$ on
    line~\ref{alg:phase1:line:initialdelaunay}. Such a tree contains
    $\eta=2^\tau=\mathcal O(2^{p^2\log(\kappa^{-1})})$ leaves.
  \end{proof}
\end{corollary}

Corollary~\ref{corollary:phase1cellcountcomplexity} tells us that the leaf count
is exponential in the parameter dimension and polynomial in the
overlap. However, if we assume that the algorithm terminates with a given finite
leaf count then the following lemma states that a linearly proportional number
of problems will have had to be solved.

\begin{lemma}
  \label{lemma:outputcomplexity}
  The iteration count $\iota$ of Algorithm~\ref{alg:phase1} is
  $\mathcal O(\eta)$.
  \begin{proof}
    A full binary tree with $\eta$ leaves has $2\eta-1$ nodes, thus at most
    $\iota=\mathcal O(\eta)$ iterations will have occured.
  \end{proof}
\end{lemma}

We have analyzed the tree complexity alone, with disregard for the complexity of
the optimization problems that need to be solved at each iteration of
Algorithm~\ref{alg:phase1}. Unlike \cite{Bemporad2006b} where convex NLPs need
to be solved at each iteration (due to their original, implicit MPC algorithm
being non-hybrid), we must solve MINLPs whose solution time is
$\mathcal O(n^\ell2^m)$ in the worst case. However, the basic assumption is that
the problem \textit{is} solvable in the first place and that on-line computation
is offloaded to an off-line solution. Therefore, we do not consider the issue of
practically solving (\ref{eq:minlp}) for a given $\theta\in\Theta$.

\begin{theorem}
  \label{theorem:onlinefdeltacomplexity}
  The on-line evaluation complexity of $f_\delta$ is $\mathcal O(p^3)$.
  \begin{proof}
    Algorithm~\ref{alg:phase1} outputs a tree with $\eta_{\mathrm{o}}$ nodes at
    the first level followed by a binary tree thereafter, where
    $\eta_{\mathrm{o}}$ is the number of simplices generated by the Delaunay
    triangulation of $\Theta$ on
    line~\ref{alg:phase1:line:initialdelaunay}. Since a simplex in $\reals^p$
    has $p+1$ facets, there are $p+1$ inequalities to evaluate in order to check
    $\theta\in\mathcal R$. Given a depth $\tau$, there are at most
    $(\eta_{\mathrm{o}}+\tau-2)(p+1)$ inequalities to evaluate. Since
    $\tau=\mathcal O(p^2)$ by Theorem~\ref{theorem:phase1treedepthcomplexity},
    the evaluation complexity of $f_\delta$ is $\mathcal O(p^3)$.
  \end{proof}
\end{theorem}

Since the evaluation complexity of $f_\delta$ is polynomial by
Theorem~\ref{theorem:onlinefdeltacomplexity}, $f_\delta$ can be used for a
guaranteed real-time warm start of an on-line mixed-integer solver.

\subsection{Choice of $\Theta$}
\label{subsec:Theta}

Throughout this paper we have assumed that $\Theta$ is available. We now explain
possible methods of obtaining it. First of all, $\Theta\subseteq\Theta^*$ should
hold. If it does not, per Theorem~\ref{theorem:phase1convergence}
Algorithm~\ref{alg:phase1} will report it in a finite number of iterations and a
different $\Theta$ must be chosen. Assuming the task of (\ref{eq:minlp}) is to
drive $\theta$ (e.g. the current state) to the origin, a $\Theta$ in a small
enough neighborhood of $0\in\reals^p$ should satisfy this property as long as
(\ref{eq:minlp}) is a well-defined controller.

Since $f_\delta$ is defined only over $\Theta$, in practice it must be ensured
that $\theta\notin\Theta$ is never encountered during runtime. This requires
$\Theta$ to be a controlled invariant set of (\ref{eq:minlp}). A straightforward
method for ensuring this is to include the constraint $\Theta^+\subseteq\Theta$
in (\ref{eq:minlp}) where $\Theta^+$ models all possible future values of
$\theta$. This is a common constraint in RMPC. In this case, which is the most
common one in practice, $\Theta$ is explicitly known. Note that convergence of
Algorithm~\ref{alg:phase1} in this case certifies recursive feasibility of
(\ref{eq:minlp}).

Finally, note from the proof of Theorem~\ref{theorem:onlinefdeltacomplexity}
that the practical complexity of $f_\delta$ is $\mathcal O(\eta_{\mathrm{o}})$
since in practice $\eta_{\mathrm{o}}\gg p$ and
$\eta_{\mathrm{o}}\gg \tau$. There is therefore a considerable interest to make
the Delaunay triangulation of $\Theta$ yield a small $\eta_{\mathrm{o}}$,
e.g. by using a simplex or a small number of simplices to define $\Theta$.

\subsection{Extensions}
\label{subsec:extensions}

The algorithm has thus far been presented as a method for partitioning
$\Theta\subseteq\Theta^*$. A possible extension is to partition the entire
$\Theta^*$. Since by Lemma~\ref{lemma:convexity} $\Theta^*_\delta$ is convex, it
may be inner-approximated with arbitrary precision \cite{Dueri2016}. Doing so
for each $\delta$, one can run the algorithm for each commutation
$\delta\in\binary^m$ by taking $\Theta=\Theta^*_\delta$. To remove overlapping
regions, the intersection of $\Theta^*_\delta$ with all previously considered
fixed-commutation feasible parameter sets can be removed.

%%% Local Variables:
%%% mode: latex
%%% TeX-master: "root"
%%% End:

\section{Illustrative Example}
\label{sec:example}

This section tests Algorithm~\ref{alg:phase1} on a set of randomly generated
dynamical systems. The goal is to demonstrate robustness by showing that the
algorithm runs successfully for a generic system and to verify the complexity
analysis of Section~\ref{subsec:complexity}.

\subsection{MPC Problem Instance Generation}

We first explain how the MPC problem instance is created for a randomly
generated dynamical system. Consider the following multiple degree-of-freedom
(DOF) oscillator, in continuous time (time is omitted for notational simplicity)
and in its configuration basis:
\begin{equation}
  \label{eq:mdof_config}
  M\ddot r+C\dot r+Kr = Lu,
\end{equation}
where $M\in\reals^{n_r\times n_r}$, $M\succ 0$, is the mass matrix,
$C\in\reals^{n_r\times n_r}$ is the damping matrix,
$K\in\reals^{n_r\times n_r}$, $K\succeq 0$, is the stiffness matrix,
$L\in\reals^{n_r\times n_u}$ is an input map, $r\in\reals^{n_r}$ is a vector of
generalized coordinates and $u\in\reals^{n_u}$ is the input. The task is to
generate $M,C,K$ and $L$ such that the system is controllable and has poles
located in a prescribed region of the complex plane. Assuming that Caughey's
condition holds \cite{SrikanthaPhani2003} such that $M,C,K$ are simultaneously
diagonalizable, (\ref{eq:mdof_config}) can be written in its modal basis:
\begin{equation}
  \label{eq:mdof_modal}
  \ddot\eta+\Lambda\dot\eta+\Omega\eta = \Gamma u,
\end{equation}
where $\eta=T^\transp M^{1/2}r\in\reals^{n_r}$ are the modal coordinates,
$T\in\reals^{n_r\times n_r}$ is the modal matrix,
$\Lambda=T^\transp M^{-1/2}CM^{-1/2}T$, $\Omega = T^\transp M^{-1/2}KM^{-1/2}T$
and $\Gamma=T^\transp M^{-1/2}L$. The matrices $\Lambda$ and $\Omega$ are
diagonal such that each row of (\ref{eq:mdof_modal}) is a 1-DOF oscillator which
contributes two poles to the overall system:
\begin{equation}
  \label{eq:3}
  \ddot\eta_i+2\zeta_i\omega_{\mathrm{n},i}\dot\eta_i+\omega_{\mathrm{n},i}^2\eta_i = u_i\quad i=1,\dots,n_r,
\end{equation}
where we chose $\Gamma=I_{n_r}$ which ensures that (\ref{eq:mdof_config}) is
controllable, $\zeta_i$ is the damping ratio and $\omega_{\mathrm{n},i}$ is the
natural frequency of the $i$-th mode. Note that this implies $n_u=n_r$. To
generate (\ref{eq:mdof_config}), it remains to choose $M$, $T$, $\Lambda$ and
$\Omega$. We choose a uniform random diagonal
$M=\mathrm{diag}(\{m_i\in [0.1,1]\}_{i=1}^{n_r})$, $T$ as the orthogonal matrix
from QR decomposition of a Gaussian random matrix and $\Lambda,\Omega$ from
randomly generating poles in the $s$-plane such that 1) the damping rate is
$\in [1,10]$ s and 2) the damped frequency is $\le 2\pi$ rad/s and 3) damping
ratio $\zeta\le 1$ with a probability of $0.2$ that $\zeta=1$.

Once generated, the system (\ref{eq:mdof_config}) is written in state space
form $\dot x = Ax+Bu+Ew$:
\begin{equation}
  \label{eq:mdof_ss}
  \hspace{-4mm}
  \setlength\arraycolsep{1pt}
  \begin{bmatrix}
    \dot r \\ \ddot r
  \end{bmatrix} =
  \begin{bmatrix}
    0 & I_{n_r} \\
    -M^{-1}K & -M^{-1}C
  \end{bmatrix}
  \begin{bmatrix}
    r \\ \dot r
  \end{bmatrix}
  +
  \begin{bmatrix}
    0 \\ M^{-1}L
  \end{bmatrix}u
  +
  \begin{bmatrix}
    0 \\
    I_{n_r}
  \end{bmatrix}
  w,
\end{equation}
where $w\in\reals^{n_r}$ is an exogenous disturbance force acting along each
generalized coordinate. This system is discretized via zero-order hold with
sampling frequency $\omega_{\mathrm{s}}=10\max_{\lambda\in\Spec(A)}|\lambda|$,
i.e. ten times faster than the fastest natural frequency present in the system
\cite{Antsaklis2007}.

Following Section~\ref{subsec:Theta}, $\Theta$ is chosen to be the smallest
robust invariant set for (\ref{eq:mdof_ss}) using the uncertainty set
$\mathcal W\triangleq\{w:\|w\|_\infty\le 10^{-3}\}$ and an LQR controller with a
$Q_{\mathrm{lqr}}=0.1I_{2n_r}$ state penalty and an $R_{\mathrm{lqr}}=I_{n_u}$
input penalty \cite{Hennet1995,Trodden2016}. For the MPC law, the uncertainty
model is changed to be norm-bounded:
\begin{equation}
  \label{eq:norm_bounded_uncertainty}
  \mathcal W'\triangleq\{w\in\reals^{n_r}:\|w\|_2\le 0.4\cdot\frac{10^{-3}\|x\|_2}{\max_{v\in\mathcal V(\Theta)}\|v\|_2}\},
\end{equation}
which is a smaller uncertainty but, importantly, introduces second-order cone
constraints into (\ref{eq:minlp}) \cite{Malyuta2019}. The ad hoc factor of $0.4$
is used to reduce uncertainty such that a planning horizon of $N=3$ is feasible
for the robust MPC law, whereas only $N=1$ is guaranteed by the computation
method for $\Theta$ \cite{Trodden2016}. Finally, to make the control problem
mixed-integer the control input is constrained to be in a non-convex set:
\begin{align}
  \label{eq:input_constraint}
  u &\in \mathcal U=\{0\}\cup(\mathcal U_{\mathrm{ext}}\setminus \mathcal U_{\mathrm{int}}), \\
  \mathcal U_{\mathrm{ext}} &\triangleq \{u\in\reals^{n_r}:-u_{\max}\le u\le u_{\max}\}, \nonumber \\
  \mathcal U_{\mathrm{int}} &\triangleq \{u\in\reals^{n_r}:-10^{-3}u_{\max}\le u\le 10^{-3}u_{\max}\}, \nonumber
\end{align}
where $u_{\max,i}=\max_{v\in\mathcal V(\Theta)}|e_i^\transp K_{\mathrm{lqr}}v|$
is the largest input magnitude required by the LQR controller along the $i$-th
generalized coordinate. Note that since $\mathcal U_{\mathrm{ext}}$ and
$\mathcal U_{\mathrm{int}}$ are origin-centered hyperrectangles, one can write
$\mathcal U=\cup_{i=1}^{2n_r+1}\mathcal U_i$ where $\mathcal U_i$ are convex
polytopes and $\mathcal U_1=\{0\}$. There are then $N(2n_r+1)$ degrees of
freedom to choose which convex subsets of $\mathcal U$ the control inputs are to
be in. The robust MPC law is then:
\begin{equation}
  \label{eq:mpc_law}
  \hspace{-4mm}
  \optimize[
    x_0 = \theta, \\
    & x_{k+1} = Ax_k+Bu_k+Ew_k\quad k=0:N-1, \\
    & x_k\in\Theta\,\,\forall w_j\in\mathcal W'\,\,\,\, j=0:k-1,k=1:N, \\
    & u_k\in\mathcal U_i\text{ for some }i\in\{1:2n_r+1\}, k=0:N-1,
  ]{}{\min}{x_k,u_k}{\sum_{k=0}^{N-1}u_k^\transp R_{\mathrm{lqr}}u_k
    +x_{k+1}^\transp Q_{\mathrm{lqr}}x_{k+1}}
  \hspace{-10mm}
\end{equation}
which can be transformed into the form (\ref{eq:minlp}) via H\"older's
inequality as was shown in \cite{Malyuta2019}. Note that $p=2n_r$, therefore
(\ref{eq:mpc_law}) is constrained to even parameter dimensions.

\subsection{Algorithm Performance Statistics}

Algorithm~\ref{alg:phase1} is applied to 100 randomly generated instances of
(\ref{eq:mpc_law}) with $n_r=1,2,3$. This demonstrates that the algorithm can
scale up to at least $p=6$ and $m=21$, with higher dimensions likely possible as
discussed in Section~\ref{sec:future}.

\begin{figure}
  \centering
  \begin{subfigure}[t]{.5\columnwidth}
    \centering
    \includegraphics[width=1\linewidth]{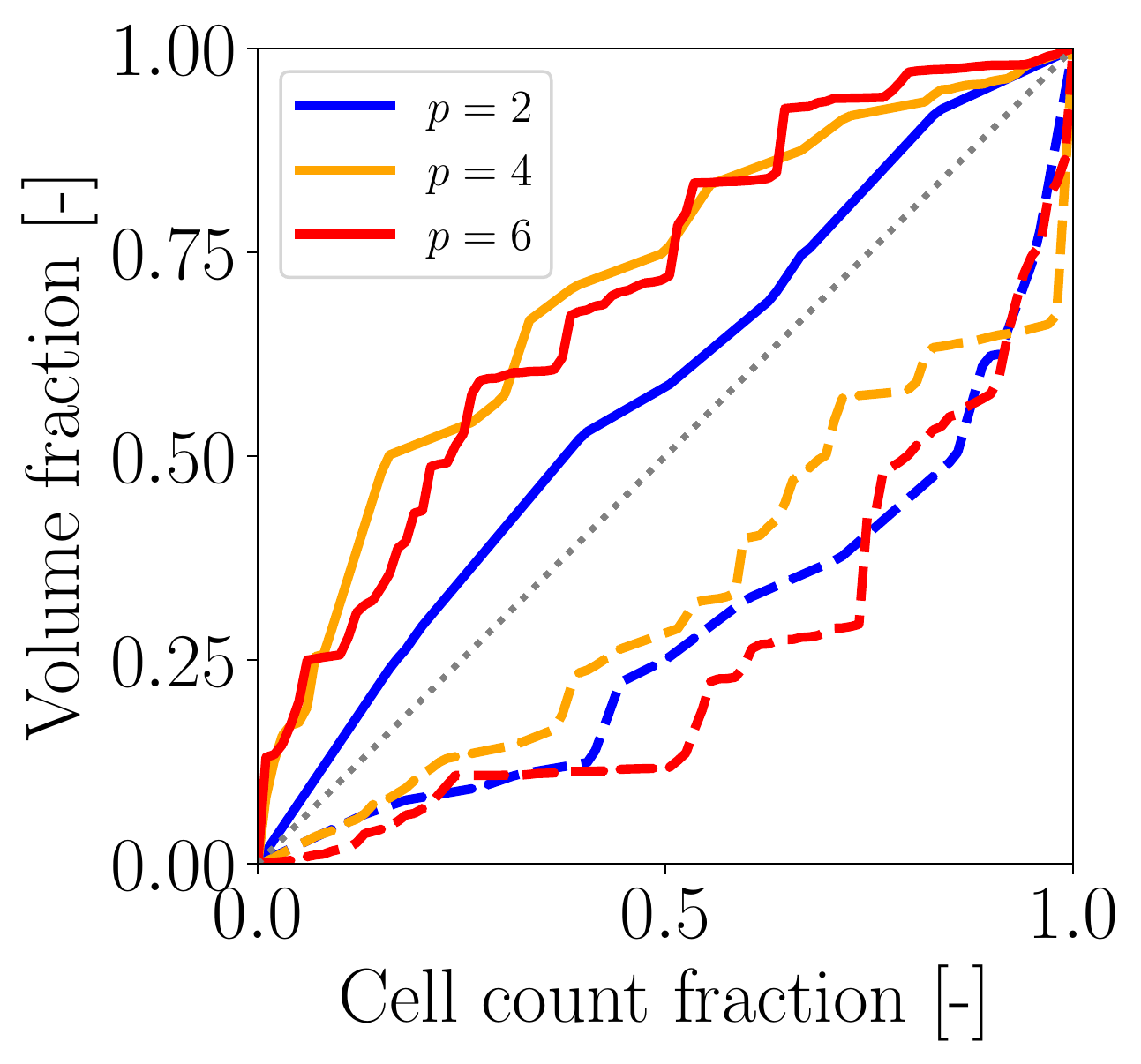}
    \caption{Volume versus cell count.}
    \label{fig:vol_vs_cell}
  \end{subfigure}%
  \begin{subfigure}[t]{.5\columnwidth}
    \centering
    \includegraphics[width=1\linewidth]{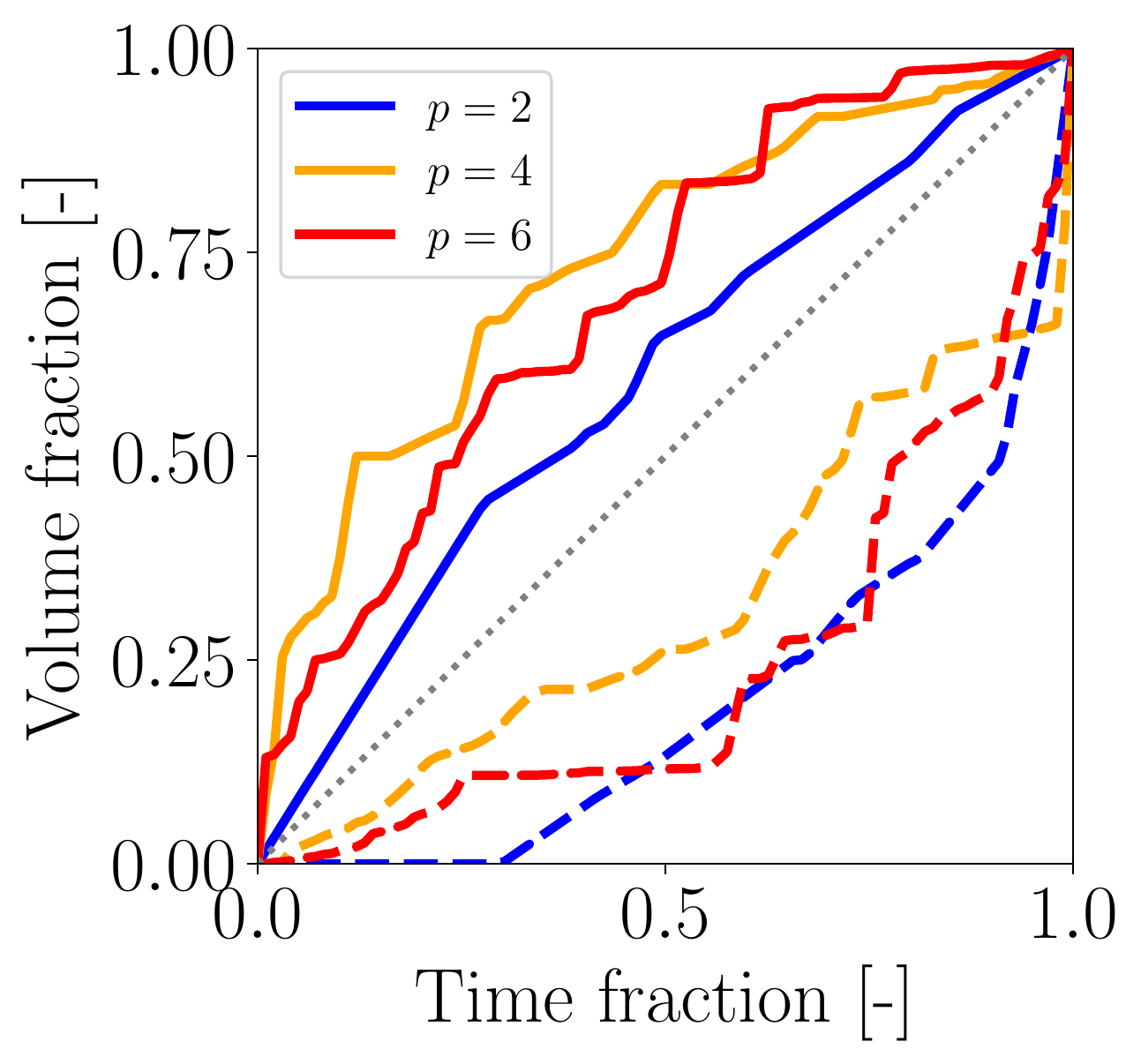}
    \caption{Volume versus runtime.}
    \label{fig:vol_vs_time}
  \end{subfigure}
  \caption{Normalized convergence plots. Cumulative closed volume, cumulative
    closed leaf count and runtime are normalized by their final respective
    values. The solid/dashed lines show the envelope max/min while the dotted
    reference line shows linear convergence.}
  \label{fig:convergence}
\end{figure}

\begin{figure}
  \centering
  \begin{subfigure}[t]{.5\columnwidth}
    \centering
    \includegraphics[width=1\linewidth]{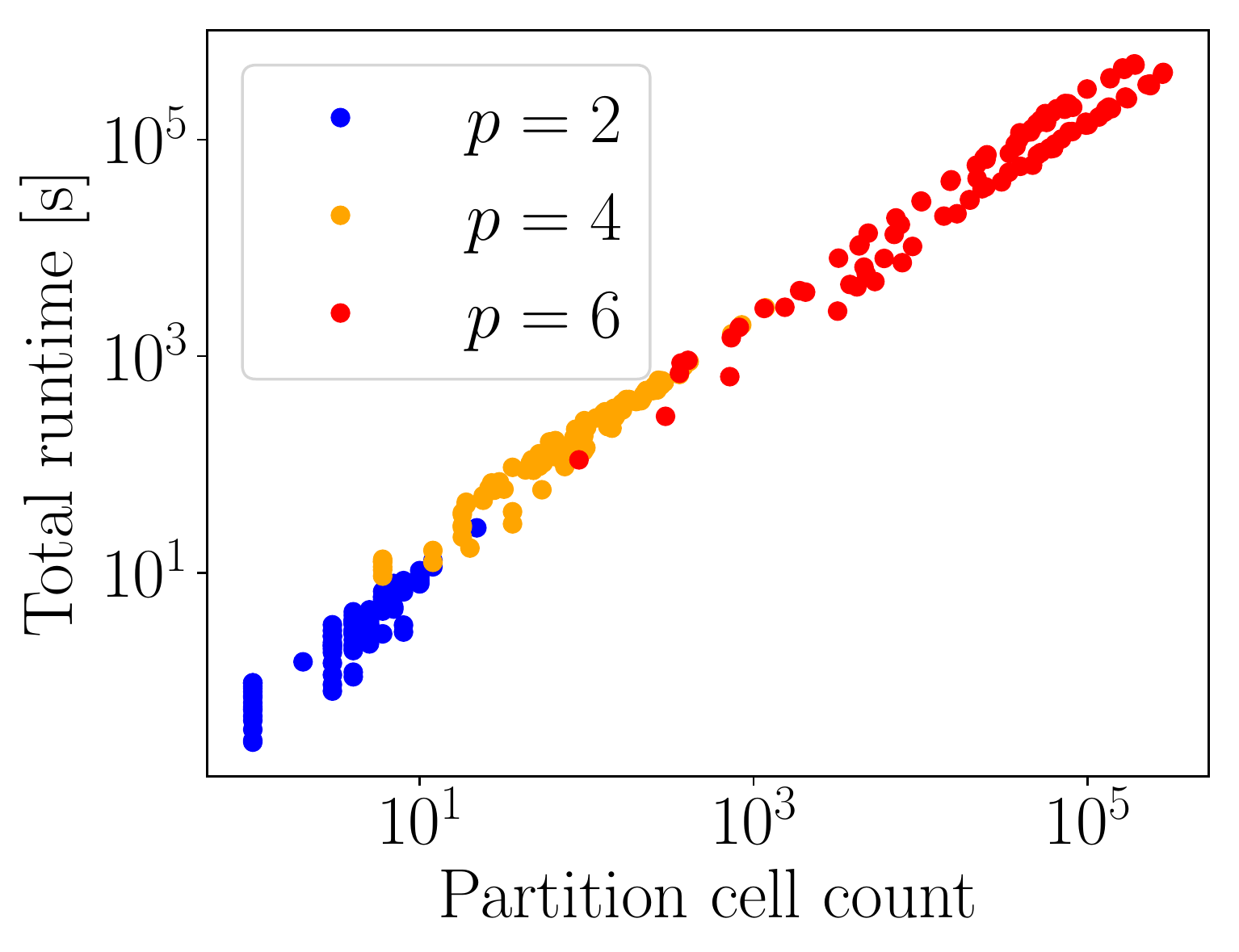}
    \caption{Total runtime statistics.}
    \label{fig:runtime_scatter}
  \end{subfigure}%
  \begin{subfigure}[t]{.5\columnwidth}
    \centering
    \includegraphics[width=1\columnwidth]{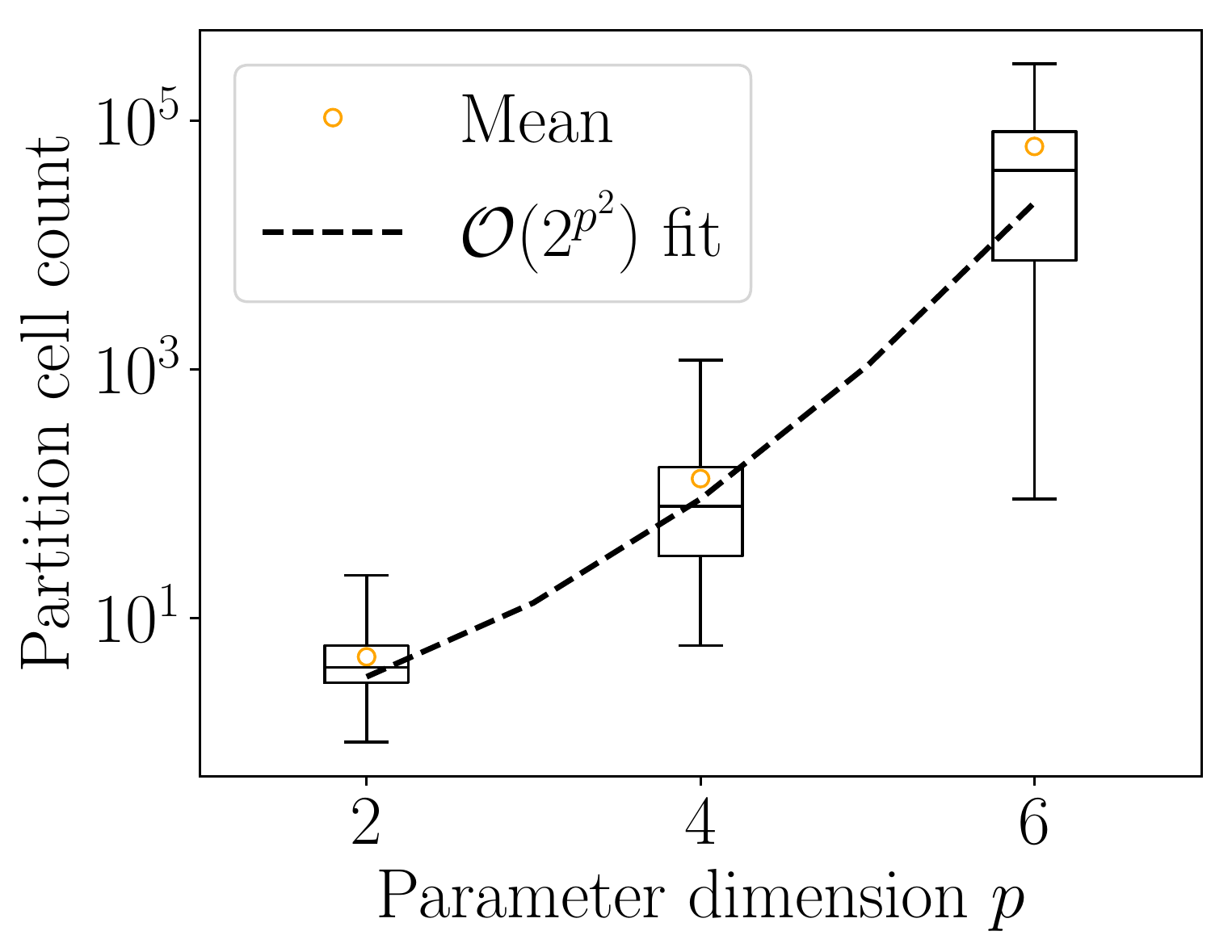}
    \caption{Cell count statistics.}
    \label{fig:cell_boxplot}
  \end{subfigure}

  \begin{subfigure}[t]{.5\columnwidth}
    \centering
    \includegraphics[width=1\linewidth]{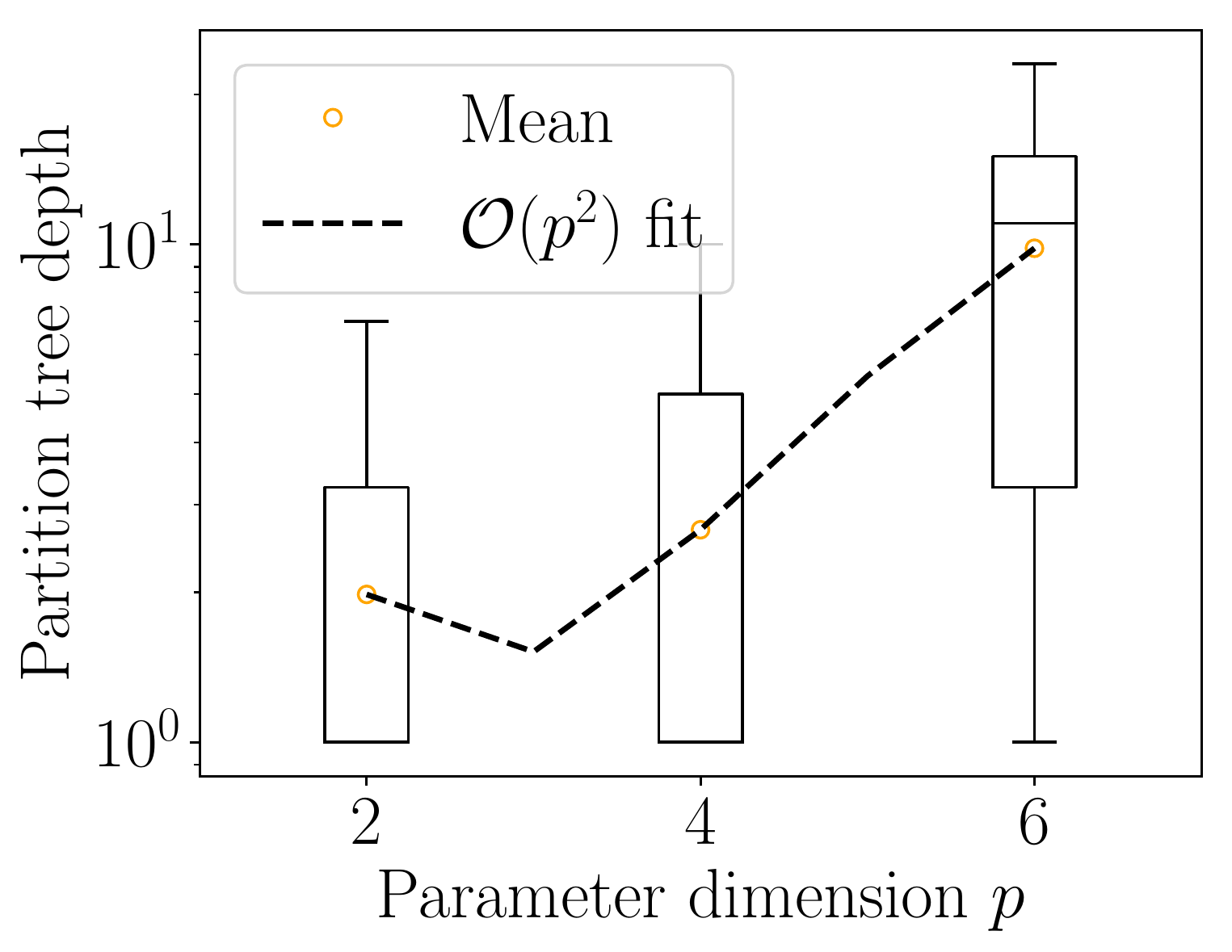}
    \caption{Tree depth statistics.}
    \label{fig:depth_boxplot}
  \end{subfigure}%
  \begin{subfigure}[t]{.5\columnwidth}
    \centering
    \includegraphics[width=1\linewidth]{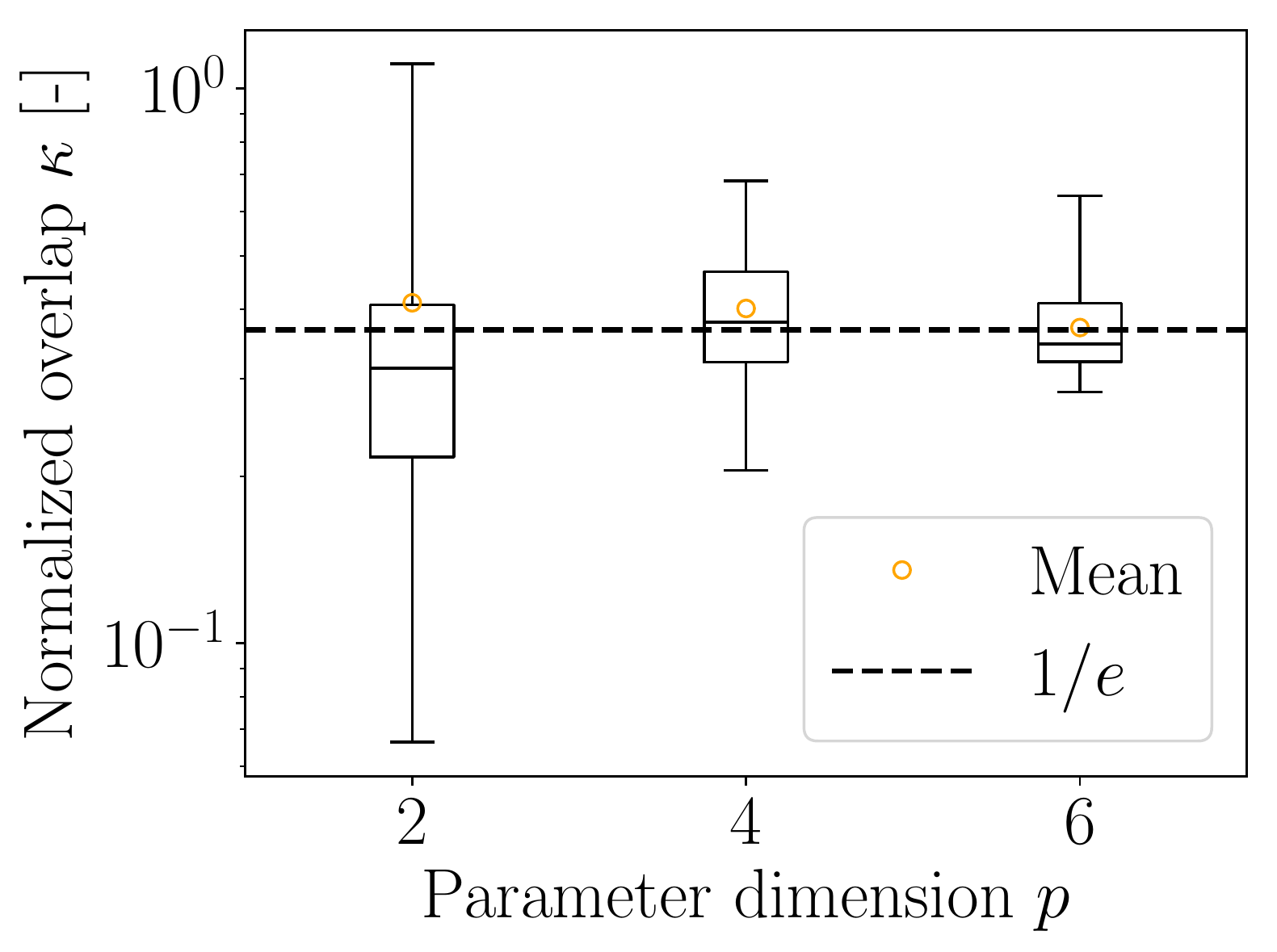}
    \caption{Normalized overlap statistics.}
    \label{fig:kappa_boxplot}
  \end{subfigure}
  \caption{Final partition statistics. Whiskers show the full range.}
  \label{fig:statistics}
\end{figure}

Figure~\ref{fig:convergence} shows convergence plots using the fraction of the
volume of $\Theta$ made up by the closed leafs as the metric. Since at first
none and in the end all leaves are closed, this metric goes from $0$ at the
start to $1$ at the end of Algorithm~\ref{alg:phase1} and is easily evaluated
since the volume of a simplex is well known \cite{Stein1966}. The algorithm has
a favorable convergence characteristic in that no convergence curve in our tests
deviated significantly from a linear rate. The practical significance of this is
that the algorithm progresses steadily towards filling up the entire volume of
$\Theta$ rather than being very slow at the beginning and very fast at the end
or vice versa. Both cases would be poor for the user to supervise since,
especially for $p>3$, it becomes difficult to diagnose the reason for slow
convergence.

Figure~\ref{fig:runtime_scatter} shows the wall-clock total runtime
corresponding to each run, which appears to increase linearly and with unity
slope as a function of the final partition cell count. The linear trend agrees
with the linear output complexity of Lemma~\ref{lemma:outputcomplexity} while
the unity slope may be interpreted as that, regardless of $p$, it takes the
current implementation on average 1 second to add 1 closed leaf to the partition
tree. Note that this measurement includes the time taken to traverse potentially
many layers of the tree until adding a closed leaf (up to about 20 layers
according to Figure~\ref{fig:depth_boxplot}). The fact that $p=2,4,6$ all lie on
the same trend line indicates that the current implementation's bottleneck is
not the complexity of the intermediate MINLPs that have to be solve but rather
e.g. database access speed.

Figures~\ref{fig:cell_boxplot}, \ref{fig:depth_boxplot} and
\ref{fig:kappa_boxplot} show statistics on the final tree leaf count and depth
along with fitted complexity curves resulting from
Section~\ref{subsec:complexity}. Figure~\ref{fig:cell_boxplot} shows clearly
that the partition tree leaf count increases exponentially with $p$ as
stipulated by
Corollary~\ref{corollary:phase1cellcountcomplexity}. Interestingly, we note from
Figure~\ref{fig:depth_boxplot} that for some instances of (\ref{eq:mpc_law}) the
tree depth does not go beyond the first layer. In other words, sometimes the
Delaunay triangulation of $\Theta$ on line~\ref{alg:phase1:line:initialdelaunay}
of Algorithm~\ref{alg:phase1} suffices. Note that because the complexities in
Theorem~\ref{theorem:phase1treedepthcomplexity} and
Corollary~\ref{corollary:phase1cellcountcomplexity} depend also on the overlap
$\kappa$, which currently cannot be computed a priori, the regressions in
Figures~\ref{fig:cell_boxplot} and \ref{fig:depth_boxplot} carry an omitted
variable bias. However, Corollary~\ref{corollary:phase1cellcountcomplexity}
allows to compute normalized values for $\kappa$ by assuming that the deviations
in Figure~\ref{fig:cell_boxplot} of the actual cell count from the fitted one
are due to $\kappa$ alone. This effectively captures the normalized variation
required from $\kappa$ in order to explain the deviation of observed results
from the regressed theoretical values. This is shown in
Figure~\ref{fig:kappa_boxplot} where $\kappa=1/e$ if the match between the
fitted and predicted cell counts is perfect. As expected, the effect of $\kappa$
diminishes for higher $p$ where the exponential complexity in $p$ dominates over
the polynomial complexity in $\kappa$.

%%% Local Variables:
%%% mode: latex
%%% TeX-master: "root"
%%% End:

\section{Future Work}
\label{sec:future}

Algorithm~\ref{alg:phase1} is subject to several potential improvements. First,
it would be interesting to compute the overlap $\kappa$ given $\Theta$ and
(\ref{eq:minlp}). This way, Algorithm~\ref{alg:phase1} could be certified to
converge a priori. Next, the partitioning process is parallelizable since
lines~\ref{alg:phase1:line:getfirstopenleaf}-\ref{alg:phase1:line:closeleaf} can
be executed in parallel for different leafs. Assuming that database
communication latency can be highly optimized to the point of being negligible,
we can say that
lines~\ref{alg:phase1:line:getfirstopenleaf}-\ref{alg:phase1:line:closeleaf} can
be made to execute entirely in parallel. Amdahl's law then predicts that
$t_{\mathrm{p}}=t_{\mathrm{s}}/n_{\text{proc}}$ where $t_{\mathrm{p}}$,
$t_{\mathrm{s}}$ and $n_{\text{proc}}$ are the parallel total runtime, serial
total runtime and number of processors used respectively. The total runtime can
thus be reduced in inverse proportion to the number of processors
available. Given that modern university facilities can typically provide access
to on the order of $10^2$ processors and that another $10^2$ factor can be
achieved by using a compiled programming language, we expect that a compiled
parallel implementation can yield a speedup of at least $10^4$. According to
Figure~\ref{fig:runtime_scatter}, this means that one could compute partitions
with $p=4$ in 0.1~s, with $p=6$ in 10~s and (extrapolating) with $p=8$ in
1000~s.

%%% Local Variables:
%%% mode: latex
%%% TeX-master: "root"
%%% End:

\section{Conclusion}
\label{sec:conclusion}

This paper presented an algorithm for generating a feasible parameter set
partition applicable to hybrid MPC problems. The algorithm consists of
systematically breaking down the feasible parameter set into smaller simplices
until these can be assigned an integer solution that is feasible everywhere in
them. Convergence in a finite number of iterations was proven with novel insight
into an overlap characteristic of the MINLP. The on-line evaluation of the
partition is polynomial time and thus can be used as a guaranteed real-time warm
start of a mixed-integer solver. Extensive testing on randomly generated systems
confirmed the complexity calculations, showed favorable convergence properties
and suggests that the algorithm is robust enough to be applied on a wide variety
of hybrid MPC problems.

%%% Local Variables:
%%% mode: latex
%%% TeX-master: "root"
%%% End:

%%%%%%%%%%%%%%%%%%%%%%%%%%%%%%%%%%%%%%%%%%%%%%%%%%%%%%%%%%%%%%%%%%%%%%%%%%%%%%%%

%%%%%%%%%%%%%%%%%%%%%%%%%%%%%%%%%%%%%%%%%%%%%%%%%%%%%%%%%%%%%%%%%%%%%%%%%%%%%%%%
\section{Acknowledgment}

Part of the research was carried out at the Jet Propulsion Laboratory,
California Institute of Technology, under a contract with the National
Aeronautics and Space Administration. Government sponsorship acknowledged. The
authors would like to extend special gratitude to Daniel P. Scharf, Jack
Aldritch and Carl Seubert for their helpful insight and discussions.

\bibliographystyle{ieeetr}
\bibliography{bibliography.bib}
%%%%%%%%%%%%%%%%%%%%%%%%%%%%%%%%%%%%%%%%%%%%%%%%%%%%%%%%%%%%%%%%%%%%%%%%%%%%%%%%

\end{document}